\documentclass[leqno,12pt]{article}
\usepackage{latexsym}
\usepackage{amsmath}
\usepackage{amssymb}
\usepackage{bbm}
\usepackage{color}
\usepackage[authoryear]{natbib}

\usepackage{graphicx}
\usepackage{epsfig}

 0

\newtheorem{thm}{Theorem}[section]
\newtheorem{prop}[thm]{Proposition}

\newtheorem{lem}[thm]{Lemma}

\makeatletter\@addtoreset{equation}{section}\makeatother

\def\proof{\noindent{\bf Proof:}\hskip10pt}
\def\QED{\hfill $ \Box $\\}

\begin{document}

\title{Nonstandard limit theorems and large deviations for the Jacobi beta ensemble}

\author{Jan Nagel}

\maketitle
\begin{abstract}
In this paper we show weak convergence of the empirical eigenvalue distribution and of the weighted spectral measure of the Jacobi ensemble, when one or both parameters grow faster than the dimension $n$. In these cases the limit measure is given by the Marchenko-Pastur law and the semicircle law, respectively. For the weighted spectral measure we also prove large deviation principles under this scaling, where the rate functions are those of the other classical ensembles.
\end{abstract}

\medskip

\section{Introduction}
The Jacobi ensemble is one of the central distributions of Hermitian matrices and can be characterised by the eigenvalue density
\begin{align}\label{jacobi}
c_{a,b,\beta} \cdot \prod_{i<j} |\lambda_i-\lambda_j|^\beta \prod_{i=1}^n \lambda_i^{a-e_\beta} (1-\lambda_i)^{b-e_\beta} \mathbbm{1}_{\{ 0<\lambda_i<1 \} }
\end{align}
with parameters $a,b>0$, $\beta>0$ and $e_\beta=\tfrac{\beta}{2}(n-1)+1$ and $c_{a,b,\beta}$ the normalisation constant. If $X,Y$ are $n\times \tfrac{2}{\beta} a$ and $n\times \tfrac{2}{\beta} b$ ($\tfrac{2}{\beta} a,\tfrac{2}{\beta} b \geq n$) matrices with i.i.d. real ($\beta=1$), complex ($\beta=2$) or quaternion ($\beta=4$) standardnormal distributed entries, then the matrix $A=(XX^* + YY^*)^{-1/2} XX^* (XX^* + YY^*)^{-1/2}$ has eigenvalues with density \eqref{jacobi}. Due to this construction, results about the spectral properties of the Jacobi ensemble can be applied in the multivariate analysis of variance (MANOVA, see \cite{muirhead1982}). In statistical mechanics, the Jacobi ensemble arises as a model for a log gas with logarithmic interactions, confined to the interval $[0,1]$ (\cite{dyson1962a}, \cite{forrester2010}). One central object in the asymptotic study of random matrix ensembles is the empirical eigenvalue distribution
\begin{align}\label{empmeas}
\hat{\mu}_n = \frac{1}{n} \sum_{i=1}^n \delta_{\lambda_i} ,
\end{align}
when $\delta_z$ denotes the Dirac measure in $z$. Recently, a number of authors showed interest in the so-called spectral measure of random matrices, a weighted version of the empirical eigenvalue distribution. For the matrix $A$ the spectral measure $\mu_n$ is defined by the relation
\begin{align} \label{moments}
\int x^k \mu_n(dx)  = \langle e_1,A^k e_1 \rangle 
\end{align}
for all $k\geq 1$, the functional calculus in the sense of the spectral theorem (\cite{dunfschw1963}). For a matrix of the classical ensembles, $\mu_n$ has the representation 
\begin{align} \label{spmeas}
\mu_n= \sum_{i=1}^n w_i \delta_{\lambda_i} ,
\end{align}
with weights $w_i=|\langle e_1,u_i\rangle|^2$, when $u_1,\dots ,u_n$ is a corresponding system of eigenvectors. Due to the invariance of the ensembles under unitary conjugations, the weights are independent from the eigenvalues and follow a Dirichlet distribution with parameter $\tfrac{\beta}{2}n$ (\cite{dawid1977}). Since both the eigenvalue and the weight distribution can be defined not only for $\beta \in \{1,2,4\}$, we will consider general $\beta >0$ and denote the joint distribution of eigenvalues \eqref{jacobi} and independent Dirichlet weights by $\mathcal{J}\beta \mathcal{E}(a,b)$.

The large $n$ behaviour of the eigenvalues with density \eqref{jacobi} has been subject to intensive research. There are numerous results such as almost sure limits (\cite{col2005}, \cite{leff1964}), CLTs (\cite{jiang2009}, \cite{nagwad1993}, \cite{johnstone2008}) or large deviations (\cite{hiaipetz2006}) about the eigenvalues in the bulk, i.e. for $\mu_n$ or of the edge behaviour of the largest eigenvalue. For a general overview we refer to the books of \cite{agz2010}, \cite{forrester2010} and \cite{mehta2004}. For the weighted random measure, large deviation principles were proven by \cite{gamrou2009a,gamrou2009b} and central limit theorems can be found in the papers of \cite{lytpast2009,lytpast2009b} and \cite{detnag2012}.

The mentioned results assume that the parameters $a=a_n$ and $b=b_n$ satisfy a \emph{standard} behaviour,
\begin{align*}
\frac{a_n}{n} \xrightarrow[n \rightarrow \infty ]{} a_0 \in[0,\infty), \quad \frac{b_n}{n} \xrightarrow[n \rightarrow \infty ]{} b_0 \in[0,\infty) .
\end{align*}
Significantly less is known in the what we call \emph{nonstandard} case, where $a_n$ and/or $b_n$ are of larger order than $n$. Under restrictions on the minimal and maximal order of parameters, \cite{detnag2009} showed that the rescaled empirical eigenvalue distribution can converge to different classes of limit distributions. In particular, these classes include the semicircle law and the Marchenko-Pastur law, the large dimensional limit of the other two classical distributions of Hermitian matrices, the Gaussian and the Laguerre ensemble, respectively. \cite{jiang2012} showed weak convergence to the Marchenko-Pastur law as well as convergence at the edge under the assumption $a_n =o(\sqrt{b_n}), n=o(\sqrt{b_n})$. Both of these results rely on uniform estimates, for the difference between eigenvalues and zeros of orthogonal polynomials (\cite{detnag2009}), and for the difference between densities (\cite{jiang2012}).\\
In this paper, we use a different approach and show weak convergence of rescaled versions of $\mu_n$ and $\hat{\mu}_n$ to the semicircle law and the Marchenko-Pastur law without any additional assumptions on the rate of the parameters. The main idea is to study the spectral measure of a tridiagonal representation of the Jacobi ensemble. The weak convergence follows from the entrywise convergence of the tridiagonal matrix, thus eliminating any need for uniform estimates. 
Furthermore, we prove large deviation principles for the spectral measure $\mu_n$ under nonstandard scaling. These LDPs illustrate further the the connection between the classical ensembles, as the rate functions are given by the rate functions of the Gaussian and the Laguerre ensemble, respectively.

This paper is structured as follows: In Section 2, we state the results and introduce the other classical ensembles. Section 3 introduces the tridiagonal models and the relation to orthogonal polynomials. Finally, Section 4 contains the proofs and some auxiliary results.

\section{Results}

To formulate the results and put them into context, we define the two remaining classical ensembles. The Gaussian (beta-)ensemble is defined by the eigenvalue density
\begin{align} \label{ewge}
f_G(\lambda) = c_G \prod_{i<j} |\lambda_i-\lambda_j|^\beta \prod_{i=1}^n e^{-\frac{1}{2}\lambda_i^2}.
\end{align}
This is the eigenvalue distribution of a Hermitian matrix with independent Gaussian entries on and above the diagonal. For growing dimension, the celebrated result of \cite{wigner1955} states that the empirical eigenvalue measure $\hat{\mu}_n$ of the rescaled eigenvalues $\sqrt{\tfrac{2}{\beta n}}\lambda_i$ converges almost surely weakly to the semicircle law $\operatorname{SC}$ with density
\begin{align}\label{semicircle}
 \frac{1}{2\pi} \sqrt{4-x^2}\mathbbm{1}_{\{ -2<x<2\} } .
\end{align}

The Laguerre (beta-)ensemble, or Wishart ensemble, has eigenvalue density
\begin{align} \label{ewle}
f_L(\lambda) = c_L^\gamma \prod_{i<j} |\lambda_i-\lambda_j|^\beta \prod_{i=1}^n \lambda_i^{a-e_\beta} e^{- \lambda_i} \mathbbm{1}_{\{ \lambda_i>0\} },
\end{align}
with parameter $a>0$, which corresponds to a square of a Gaussian matrix. As $n\to \infty$ and $\tfrac{\beta n}{2a}\to \tau \in (0,1]$, the measure $\hat{\mu}_n$ of the rescaled eigenvalues $\tilde{\lambda}_i= \tfrac{2\tau}{n\beta}\lambda_i$ tends almost surely to the Marchenko-Pastur law $\operatorname{MP}(\tau)$ with parameter $\tau$ and density
\begin{align}\label{marchpast}
 \frac{1}{2\pi \tau x}\sqrt{(x-\tau^-)(\tau^+-x)}  \mathbbm{1}_{\{ \tau^-<x<\tau^+\} },
\end{align}
where $\tau^- = (\sqrt{\tau}-1)^2$ and $\tau^+ = (\sqrt{\tau}+1)^2$.

\subsection{Weak convergence}

Our first two main theorems deal with the almost sure weak convergence of the rescaled eigenvalue distribution of the Jacobi ensemble. In this case, the limit is the same for the empirical measure $\hat\mu_n$ and the spectral measure $\mu_n$. In what follows, we regard the random measures as elements of the set $\mathcal{M}_c$ of probability measures on $\mathbb{R}$ with compact support, endowed with the weak topology and the corresponding Borel $\sigma$-algebra. To simplify notation, we write $\beta'=\tfrac{1}{2} \beta$.

\begin{thm}\label{LLN2}
Suppose in the Jacobi ensemble $\mathcal{J}\beta\mathcal{E}(a_n,b_n)$ the parameters satisfy
\begin{align*}
\frac{a_n}{n} \xrightarrow[n \rightarrow \infty ]{} \infty,\quad \frac{a_n}{b_n} \xrightarrow[n \rightarrow \infty ]{} \sigma >0 
\end{align*}
and additionally
\begin{align*}
\sqrt{\frac{b_n}{n}}\left(\sigma-\frac{a_n}{b_n}\right) \xrightarrow[n \rightarrow \infty ]{} 0.
\end{align*}
Then the spectral measure $\mu_n$ and the empirical eigenvalue distribution $\hat{\mu}_n$ of the rescaled eigenvalues 
\begin{align*}
(\sigma +1) \sqrt{b_n\frac{1+\sigma}{\sigma  n \beta' }} \left( \lambda_i - \frac{\sigma}{\sigma+1} \right)
\end{align*} 
converge weakly to the semicircle law $\operatorname{SC}$ almost surely.
\end{thm}

\medskip

\begin{thm} \label{LLN1}
If the parameters of the Jacobi ensemble $\mathcal{J}\beta\mathcal{E}(a_n,b_n)$ satisfy 
\begin{align*}
\frac{\beta' n}{a_n} \xrightarrow[n \rightarrow \infty ]{} \tau \in (0,1],\quad \frac{b_n}{n} \xrightarrow[n \rightarrow \infty ]{} \infty ,
\end{align*}
then the spectral measure $\mu_n$ and the empirical eigenvalue distribution $\hat{\mu}_n$ of the rescaled eigenvalues 
\begin{align*}
\frac{b_n}{a_n} \lambda_i
\end{align*} 
converge weakly to the Marchenko-Pastur law $\operatorname{MP}(\tau)$ almost surely.
\end{thm}

\medskip

The fact that the rescaled Jacobi eigenvalue distribution can have the same large dimensional limit as the Gaussian and the Laguerre ensemble is closely related to the finite dimensional convergence to the other ensembles. To explain this connection, we state here the convergence for fixed dimension. These weak convergence results for fixed size are probably well-known (see e.g. \cite{gnrw2012} for $\beta=2$), however, we could not find a general reference for them. 

\begin{prop} \label{findim} 
\begin{itemize}
\item[(i)] Let $\lambda_N$ be an eigenvalue vector of the Jacobi ensemble $\mathcal{J}\beta\mathcal{E}(a_N,b_N)$ with density \eqref{jacobi}. Suppose $a_N\to \infty$ and $a_N/b_N \to \sigma \in (0,\infty)$ while the dimension $n$ is fixed and additionally
\begin{align*}
\sqrt{a_N}\left(\sigma-\frac{a_N}{b_N}\right) \xrightarrow[N \rightarrow \infty ]{} 0,
\end{align*}
then, with $1_n=(1,\dots ,1)$, 
\begin{align*}
(\sigma +1) \sqrt{b_N \frac{\sigma+1}{\sigma}} \left( \lambda_N - \frac{\sigma}{\sigma+1} 1_n \right)
\end{align*}
converges in distribution to a random vector with density $\eqref{ewge}$.
\item[(ii)] Let $\lambda_N$ be an eigenvalue vector of the $\mathcal{J}\beta\mathcal{E}(a,b_N)$ with density \eqref{jacobi}. If $b_N\to \infty$ with fixed dimension $n$, then the rescaled eigenvalue vector 
\begin{align*}
b_N\lambda_N
\end{align*}
converges in distribution to a random vector with density $\eqref{ewle}$.
\end{itemize}
\end{prop}

\medskip
Heuristically, Theorem \ref{LLN1} is now a combination of part $(ii)$ of the Proposition with the large dimensional limit of the Laguerre ensemble. We let $b_n$ tend to infinity faster than $n$ to approach the Laguerre ensemble, while simultaneously increasing dimension and choosing $a_n$ such that the Laguerre eigenvalues converge to the right limit. Theorem \ref{LLN2} has a similar explanation: parameters $a_n$ and $b_n$ of the same order growing faster than the dimension carry part $(i)$ of the proposition over in the large dimensional limit.

\medskip

\subsection{Large deviations}

We first recall the definition of a large deviation principle. Let $U$ be a topological Hausdorff space with Borel $\sigma$-algebra $\mathcal{B}(U)$. We say that a sequence $(P_{n})_n$ of probability measures on $(U,\mathcal{B}(U))$ satisfies a large deviation principle (LDP) with speed $a_n$ and
rate function $\mathcal{I}$  if:
\begin{itemize}
\item[(i)] For all closed sets $F \subset U$:
\begin{align*}
\limsup_{n\rightarrow\infty} \frac{1}{a_n} \log P_{n}(F)\leq -\inf_{x\in F}\mathcal{I}(x)
\end{align*}
\item[(ii)] For all open sets $O \subset U$:
\begin{align*}
\liminf_{n\rightarrow\infty} \frac{1}{a_n} \log P_{n}(O)\geq -\inf_{x\in O}\mathcal{I}(x)
\end{align*}
\end{itemize}
The rate function $\mathcal{I}$ is good if its level sets
$\{x\in U |\ \mathcal{I}(x)\leq a\}$ are compact for all $a\geq 0$. In our case, the measures $P_n$ will be the distributions of the random spectral measures $\mu_n$ and we will say that the sequence of measures $\mu_n$ satisfies an LDP. Similar to the weak convergence results in Section 2.1, the rate function is determined by the large deviation behaviour of the two other classical ensembles. To formulate the rate functions we first need some definitions. The Kullback-Leibler distance between two probability measures $\mu$ and $\nu$ is given by
\begin{align*}
\mathcal{K}(\mu|\nu) = \int \log \left( \frac{\partial \mu}{\partial \nu} \right) d\mu 
\end{align*}
if $\mu$ is absolutely continuous with respect to $\nu$ and $\mathcal{K}(\mu|\nu)=\infty$ otherwise. Let for $x\geq 2$
\begin{align*}
\mathcal{F}_G(x) = \int_2^x \sqrt{t^2-4} dt = \tfrac{x}{2} \sqrt{x^2-4} - 2 \log \left( \tfrac{x+\sqrt{x^2-4}}{2}\right)
\end{align*}
and define $\mathcal{F}_G(x)=\mathcal{F}_G(-x)$ for $x\leq -2$. Note that $\mathcal{F}_G$ is the rate function for the largest eigenvalue of the Gaussian ensemble (see \cite{agz2010}). Following \cite{simon05}, we say that a probability measure $\mu$ satisfies the Blumenthal-Weyl condition (B.W.c.) if
\begin{itemize}
\item[(i)] $\operatorname{supp}(\mu) = [-2,2] \cup \{E_j^-\}_{j=1}^{N^-} \cup \{E_j^+\}_{j=1}^{N^+}$, where $N^-,N^+$ may be 0, finite or infinite with
\begin{align*}
E_1^-<E_2^-<\dots <-2 \quad \text{and} \quad E_1^+>E_2^+>\dots >2 .
\end{align*}
\item[(ii)] If $N^-$ or $N^+$ is infinite, then $E_j^-$ converges towards $-2$ and $E_j^+$ converges to 2, respectively.
\end{itemize}
\cite{gamrou2009b} showed that the sequence of weighted spectral measures $\mu_n $ of the Gaussian ensemble satisfies the LDP with speed $\beta'n$ and good rate function
\begin{align} \label{rateG}
\mathcal{I}_G(\mu) = \mathcal{K}(\operatorname{SC}|\mu) + \sum_{j=1}^{N^-} \mathcal{F}_G(E_j^-) + \sum_{j=1}^{N^+} \mathcal{F}_G(E_j^+)
\end{align}
if $\mu$ satisfies B.W.c. and $\mathcal{I}_G(\mu)=\infty$ otherwise. Using the tridiagonal representation of spectral measures we introduce in Section 3, the rate function $\mathcal{I}_G$ can be written in terms of recursion coefficients of the measure $\mu$. This large deviation behaviour is fundamentally different from the one of the empirical measure $\hat\mu_n$, which satisfies a LDP with speed $n^2$ and rate function related to Voiculescu's entropy. Our first nonstandard LDP shows the approximation of the Gaussian ensemble.

\medskip

\begin{thm} \label{LDP2}
Suppose that the parameters of the Jacobi ensemble $\mathcal{J}\beta\mathcal{E}(a_n,b_n)$ satisfy
\begin{align*}
\frac{a_n}{ n} \xrightarrow[n \rightarrow \infty ]{} \infty,\quad \frac{b_n}{n} \xrightarrow[n \rightarrow \infty ]{}\infty ,\quad \text{ and } \quad
\frac{a_n-b_n}{\sqrt{b_nn}} \xrightarrow[n \rightarrow \infty ]{} 0 .
\end{align*}
Then the spectral measure $\mu_n$ of the rescaled eigenvalues 
\begin{align*}
4 \sqrt{\frac{b_n}{ n \beta}} \left( \lambda_i - \frac{1}{2} \right)
\end{align*} 
satisfies the LDP with speed $\beta'n$ and good rate function $\mathcal{I}_G$.
\end{thm}

\medskip

In the case of the Laguerre ensemble, the sequence of spectral measures $\mu_n$ satisfies an LDP with speed $\beta'n$ and a good rate function $\mathcal{I}_L$, which is conjectured to be the Laguerre-analogue of \eqref{rateG}. However, there is an explicit formulation of $\mathcal{I}_L$ in terms of the recursion variables introduced in the following section and we refer to equation \eqref{rateL} for the definition. 

\medskip

\begin{thm} \label{LDP1}
Suppose that the parameters of the Jacobi ensemble $\mathcal{J}\beta\mathcal{E}(a_n,b_n)$ satisfy
\begin{align*}
\frac{\beta' n}{a_n} \xrightarrow[n \rightarrow \infty ]{} \tau \in (0,1],\quad \frac{b_n}{n} \xrightarrow[n \rightarrow \infty ]{} \infty ,
\end{align*}
then the spectral measure $\mu_n$ of the rescaled eigenvalues 
\begin{align*}
\frac{b_n}{a_n} \lambda_i
\end{align*} 
satisfies the LDP with speed $a_n$ and good rate function $\mathcal{I}_L$.
\end{thm}

\medskip

Note that the LDP in Theorem \ref{LDP1} for the spectral measure of the rescaled eigenvalues implies the almost sure convergence to the Marchenko-Pastur law. Thus, Theorem \ref{LLN1} is a direct consequence of Theorem \ref{LDP1} and Lemma \ref{empspectral}. The LDP in Theorem \ref{LDP2} is however more restrictive concerning the behaviour of the parameters than our convergence result in Theorem \ref{LLN2}. More precisely, the assumption $\frac{a_n-b_n}{\sqrt{b_nn}} \to 0$ implies $\sigma=1$ in Theorem \ref{LLN2}.

We remark that although the two theorems suggest an exponentially fast approximation of the Gaussian and Laguerre ensemble, this would require constraints on the order of $a_n$ and $b_n$. To completely avoid these constraints, we use a different method of proof. 

\medskip

\section{Tridiagonal Representations}

The proofs are based on tridiagonal matrix models for the classical ensembles, valid for all values of $\beta>0$ and whose entries can be decomposed into independent random variables. In the following, we say that a random variable is $\operatorname{Gamma}(a)$ distributed, if it has the density 
\begin{align*}
\frac{x^{a-1}}{\Gamma(a)} e^{-x} \mathbbm{1}_{\{ x>0\} }
\end{align*}
with $a>0$. The logarithmic moment generating function of the $\operatorname{Gamma}(a)$ distribution is $\Lambda(t)=-a \log(1-t)$ (for $t<1$) with Fenchel-Legendre transform $\Lambda^*(x)=a\cdot g(a^{-1}x)$, where the function $g$ is defined by 
\begin{align} \label{functiong}
g(x)=x-\log x-1
\end{align}
if $x>0$ and $g(x)=\infty$ otherwise.

For the Jacobi ensemble, the tridiagonal matrix model was found by \cite{kilnen2004} and is constructed as follows.
Let $p_1,\dots ,p_{2n-1}$ be independent random variables distributed as
\begin{align} \label{defbetaentries}
p_k \sim \begin{cases} \operatorname{Beta}\left( \tfrac{2n-k}{2} \beta' , a+b -\tfrac{2n+k-2}{2} \beta'  \right)\quad & k \mbox{ even,} \\
				\operatorname{Beta}\left( a-\tfrac{k-1}{2} \beta' , b- \tfrac{k-1}{2} \beta'  \right)\quad & k \mbox{ odd.} \end{cases} 
\end{align}
and define 
\begin{align} \label{zerl1}
\begin{split}
d_k =& p_{2k-2} (1-p_{2k-3}) + p_{2k-1}(1-p_{2k-2})\\
c_k =& \sqrt{ p_{2k-1}(1-p_{2k-2})p_{2k}(1-p_{2k-1})}
\end{split}
\end{align}
with $p_{-1}=p_0=0$. Then the tridiagonal matrix
\begin{align} \label{tridiagonal}
\mathcal{J}_n(\beta, a,b) = \begin{pmatrix} 
  d_1 & c_1    &         &         \\
                c_1 & d_2    & \ddots  &         \\
                    & \ddots & \ddots  & c_{n-1} \\
                    &        & c_{n-1} & d_n
\end{pmatrix} ,
\end{align}
is a matrix of the Jacobi ensemble, that is, the eigenvalues follow the density \eqref{jacobi} and the square moduli of the top entries of the eigenvectors are Dirichlet distributed with parameter $\beta'$ independent of the eigenvalues. Note that our parametrisation and scaling differs from the one used by \cite{kilnen2004}, who consider eigenvalues in $[-2,2]$ and beta distributions obtained from the one in \eqref{defbetaentries} by the transformation $x\mapsto 1-2x$.

Such a sparse matrix model is particularly convenient when we consider weak convergence of spectral measures. The $k$-th moment of the spectral measure $\mu_n$ is given by the upper left entry of the $k$-th power of the matrix, compare \eqref{moments}. Thus, it depends only on a finite number of entries, and entrywise convergence of the tridiagonal matrix to a (possible infinite) matrix implies convergence of moments and then weak convergence of the spectral measure (given that the limit measure is uniquely determined by its moments). To prove convergence of moments of the empirical eigenvalue measure, we would need to control a growing number of entries and uniform statements.

The tridiagonal representation of the two other classical ensembles were proven by \cite{dumede2002}. The tridiagonal matrix $\mathcal{G}_n(\beta)$ of the Gaussian ensemble with notation as in \eqref{tridiagonal} has standardnormal distributed diagonal entries and positive subdiagonal entries with $c_k^2\sim \operatorname{Gamma}((n-k)\beta')$, such that $d_1,\dots,d_n,c_1,\dots c_{n-1}$ are independent. As $n\to \infty$, the standardised matrix $\frac{1}{\sqrt{\beta' n}}\mathcal{G}_n(\beta)$ converges entrywise almost surely to the infinite tridiagonal matrix
\begin{align*}
\mathcal{SC} = \begin{pmatrix} 
  0 \ & 1\    &         &     \ \    \\
                1 & 0    & 1  &         \\
                    & 1 & \ddots  &\ddots  \\
                    &       &  \ddots & 
\end{pmatrix} ,
\end{align*}
which has the semicircle distribution $\operatorname{SC}$ as spectral measure.

For the matrix $\mathcal{L}_n(\beta,a)$ of the Laguerre ensemble, the entries can be written as
\begin{align} \label{zerl2}
\begin{split}
d_k =& z_{2k-2} + z_{2k-1}, \\
c_k =& \sqrt{ z_{2k-1}z_{2k}},
\end{split}
\end{align}
where $z_{0}=0$ and $z_1,\dots ,z_{2n-1}$ are independent random variables with distributions
\begin{align*}
z_{2k-1} \sim &\ \operatorname{Gamma}(a- ( k-1)\beta') , \\
 z_{2k} \sim &\ \operatorname{Gamma}( ( n-k)\beta') .
\end{align*}
Under the classical rescaling as in Section 2, the limit matrix is the tridiagonal matrix $\mathcal{MP}(\tau)$ with diagonal entries $d_1=1$, $d_k=1+\tau$ for $k>1$ and subdiagonal entries $c_k = \sqrt{\tau}$. The corresponding spectral measure is the Marchenko-Pastur law $\operatorname{MP}(\tau)$.

The weak convergence in Proposition \ref{findim} is a direct consequence of the tridiagonal representation of the ensembles and the weak convergence of the scalar beta distribution: if $X_N\sim Beta(a,b_N)$, then $b_NX_N$ converges in distribution to a $\operatorname{Gamma}(a)$ distributed random variable and if $X_N \sim \operatorname{Beta}(a_N,b_N)$ with parameters satisfying the assumptions as in case $(i)$, then
\begin{align*}
(\sigma +1) \sqrt{b_N \frac{\sigma+1}{\sigma}} \left( X_N - \frac{\sigma}{\sigma+1} \right) 
\xrightarrow[n \rightarrow \infty ]{d} \mathcal{N}(0,1) .
\end{align*}
The componentwise convergence of the tridiagonal models gives then the weak convergence of the eigenvalues.

Besides parametrizing the classical ensembles by independent random variables, the tridiagonal models have another fundamental property: they represent the multiplication $f(x)\mapsto xf(x)$ in $L^2(\mu_n)$, when $\mu_n$ is the spectral measure and the basis of $L^2(\mu_n)$ is $\{P_0,\dots P_{n-1}\}$ with $P_j$ the $j$-th orthonormal polynomial. As a consequence, the entries in \eqref{tridiagonal} appear in the three term recursion of the orthogonal polynomials,
\begin{align} \label{polrek}
xP_j(x) = c_{j+2}P_{j+1}(x) + d_{j+1}P_j(x) + c_{j} P_{j-1}(x)
\end{align}
for $1\leq j\leq n-3$. For more on this relation, we refer to \cite{simon97, simon05b}. The measure $\mu$ is supported on $[0,\infty)$ if and only if there are nonnegative $z_1,z_2,\dots$ such that the decomposition in \eqref{zerl2} holds. Furthermore, if $\mu$ is concentrated on $[0,1]$, then $z_k=p_k(1-p_{k-1})$ form a chain sequence with $p_k\in [0,1]$. In this case we get a decomposition of recursion coefficients as in \eqref{zerl1}. It is possible to formulate the rate function $\mathcal{I}_G$ of the Gaussian ensemble defined in \eqref{rateG} in terms of the recursion coefficients. The sum-rule by \cite{kilsim03} states that
\begin{align}\label{rateG2}
\mathcal{I}_G(\mu)= \sum_{k=1}^\infty \frac{1}{2} d_k^2 + g(c_k^2) ,
\end{align}
where $d_k$ and $c_k$ are the recursion coefficients of polynomials orthonormal with respect to $\mu$, $g$ is as in \eqref{functiong} and both sides may be equal to $+\infty$ simultaneously. With the decomposition of recursion coefficients, we are able to formulate the rate function $\mathcal{I}_L$ in the LDP of the spectral measure of the Laguerre ensemble. In this case, all measures are supported by $[0,\infty)$ and we have (\cite{gamrou2009b})
\begin{align}\label{rateL}
\mathcal{I}_L(\mu) =   \sum_{k=1}^\infty g(z_{2k-1}) + \tau g(z_{2k}/\tau) ,
\end{align}
with $z_k$ as in \eqref{zerl2}. If the support of $\mu$ is not a subset of $[0,\infty)$, we set $\mathcal{I}_L(\mu)=\infty$. Note that the minimum of $\mathcal{I}_L$ is attained for $z_{2k-1}=1,z_{2k}=\tau$, which corresponds to the Marchenko-Pastur law $\operatorname{MP}(\tau)$.

\medskip

\section{Proofs}

Before we start to prove the main theorems of Section 2, we first prove some auxiliary results for the concentration and large deviations of scalar beta and gamma distributed random variables.

\medskip

\begin{lem} \label{betaab} Let $X\sim \operatorname{Beta}(a,b)$ be a Beta-distributed random variable, then for any $\tfrac{1}{2}>\varepsilon>0$,
\begin{align*}
P\left( |X-E[X]| > \varepsilon \right) \leq 4 \exp \left\{ -\frac{\varepsilon^2}{128} \cdot \frac{a^3+b^3}{ab} \right\}
\end{align*}
\end{lem}

\medskip

\proof
Let $Y\sim \operatorname{Gamma}(a), Z\sim \operatorname{Gamma}(b)$ be independent gamma distributed random variables with mean $a$ and $b$, respectively. We will make use of the equality in distribution
\begin{align} \label{betafraction}
X \stackrel{d}{=} \frac{Y}{Y+Z} .
\end{align} 
The standard Chernoff bounds for gamma distributed random variables give for $\varepsilon>0$ the inequalities
\begin{align*}
P\big( Y>a(1+\varepsilon)\big)\leq \exp \big( a(\log (1+\varepsilon)-\varepsilon) \big), \qquad P\big( Y<a(1-\varepsilon)\big)\leq \exp \big( a(\log (1-\varepsilon)+\varepsilon) \big) .
\end{align*}
Using that $\log(1-\varepsilon)+\varepsilon \leq \log(1+\varepsilon)-\varepsilon\leq - \tfrac{1}{4}\varepsilon^2$ for $0<\varepsilon<\tfrac{1}{2}$, this yields
\begin{align*}
P\left( \left| \tfrac{1}{a} Y -1\right| >\varepsilon \right) \leq 2 \exp \left( -a \tfrac{\varepsilon^2}{4} \right)
\end{align*}
To get an inequality for the beta distribution note that $|A-1|<\varepsilon$ and $|B-1|<\varepsilon$ implies for any $c\geq 0$
\begin{align*}
\left| \frac{A}{cA+B} - \frac{1}{1+c} \right| \leq \frac{|A-1|+|B-1|}{(cA+B)(c+1)}\leq \frac{2\varepsilon}{(1-\varepsilon)(c+1)^2}  \leq 4 \varepsilon ,
\end{align*}
where we used $\varepsilon<\tfrac{1}{2}$ for the last inequality. Combining this implication with $c=\tfrac{a}{b}$ with the Chernoff bounds for $\tfrac{1}{a}Y$ and $\tfrac{1}{b}Z$ gives
\begin{align*}
P\left( \tfrac{b}{a}|X-E[X]| > \varepsilon \right) = P\left( \left| \frac{\frac{1}{a}Y}{ \frac{1}{b}Y+\frac{1}{b}Z} - \frac{1}{\frac{a}{b} +1} \right| > \varepsilon \right) 
\leq  2 \exp \left( -a \tfrac{\varepsilon^2}{64} \right) + 2 \exp \left( -b \tfrac{\varepsilon^2}{64} \right) .
\end{align*}
On the other hand, the symmetry of the beta distribution yields also the inequality
\begin{align*}
P\left( \tfrac{a}{b}|X-E[X]| > \varepsilon \right)  = &P\left( \tfrac{a}{b}|(1-X)-E[(1-X)]| > \varepsilon \right) \\
= & P\left( \left| \frac{\frac{1}{b}Z}{ \frac{1}{b}Y+\frac{1}{b}Z} - \frac{1}{\frac{b}{a} +1} \right| > \varepsilon \right)  \\
\leq & 2 \exp \left( -a \tfrac{\varepsilon^2}{64} \right) + 2 \exp \left( -b \tfrac{\varepsilon^2}{64} \right) .
\end{align*}
These two inequalities together give
\begin{align*}
P\left( \max \{ \tfrac{a}{b}, \tfrac{b}{a} \} |X-E[X]| > \varepsilon \right) \leq 4 \exp \left( -\min \{ a,b \} \tfrac{\varepsilon^2}{64} \right) ,
\end{align*}
which implies 
\begin{align*}
& P\left(  |X-E[X]| > \varepsilon \right) \leq 4 \exp \left( -\tfrac{\varepsilon^2}{64} \min \{ a,b \} \cdot \left( \max \{ \tfrac{a}{b}, \tfrac{b}{a} \} \right)^2  \right) \\ = 
& 4 \exp \left( -\tfrac{\varepsilon^2}{64} \max \{  \tfrac{a^2}{b}, \tfrac{b^2}{a} \}   \right)
\leq 4  \exp \left( -\tfrac{\varepsilon^2}{64} \cdot \tfrac{1}{2} \left( \tfrac{a^2}{b}+ \tfrac{b^2}{a}\right)    \right) ,
\end{align*}
and we obtain the stated inequality.
\QED

\medskip

The following theorem shows that the distance between the empirical eigenvalue distribution and the corresponding spectral measure goes to zero almost surely. Here, we consider the Kolmogorov (or uniform) distance defined for measures $\mu$, $\nu$ on the real line with distribution functions $F_\mu$ and $F_\nu$, respectively, by
\begin{align*}
d_K(\mu,\nu) = \sup_{x\in \mathbb{R}} |F_\mu(x) - F_\nu(x)| .
\end{align*}
Note that convergence in the Kolmogorov metric implies weak convergence.

\begin{thm}\label{empspectral}
For the empirical eigenvalue distribution $\hat\mu_n$ defined in \eqref{empmeas} and the corresponding spectral measure $\mu_n$ as in \eqref{spmeas}, we have
\begin{align*}
d_K(\hat\mu_n, \mu_n) \xrightarrow[n \rightarrow \infty ]{} 0
\end{align*}
almost surely.
\end{thm}

\medskip

\proof
The measures $\hat\mu_n$ and $\mu_n$ have the same support points $\lambda_1,\dots, \lambda_n$, which implies for the Kolmogorov distance
\begin{align*}
d_K(\hat\mu_n, \mu_n) = \max_{1\leq k\leq n} \left|\sum_{i=1}^k \left(w_i-\frac{1}{n}\right) \right| .
\end{align*}
Recall that the weight vector $(w_1,\dots,w_n)$ is Dirichlet distributed with parameter $\beta'$ and has the same distribution as
\begin{align*}
\left( \frac{G_1}{G_1+\dots +G_n}, \frac{G_2}{G_1+\dots +G_n}, \dots ,\frac{G_1}{G_1+\dots +G_n} \right)
\end{align*}
where $G_1,\dots ,G_n$ are independent $\operatorname{Gamma}(\beta')$ distributed. Therefore, the sum $w_1+\dots +w_k$ is $\operatorname{Beta}( k \beta',(n-k)\beta')$ distributed with mean $\tfrac{k}{n}$. Applying Lemma \ref{betaab} gives for $0<\varepsilon<\tfrac{1}{4}$
\begin{align*}
P\left( \max_{1\leq k\leq n} \left|\sum_{i=1}^k \left(w_i-\frac{1}{n}\right) \right| >\varepsilon \right) 
\leq & \sum_{k=1}^n P\left(\left| w_1+\dots +w_k-\tfrac{k}{n} \right| >\varepsilon \right) \\
\leq & \sum_{k=1}^n 4 \exp\left( -\frac{\varepsilon^2}{128} \cdot \beta' \cdot \frac{k^3+(n-k)^3}{k(n-k)} \right) \\
=& \sum_{k=1}^n 4 \exp\left( -\frac{\varepsilon^2}{128} \cdot \beta' n \cdot \frac{(\frac{k}{n})^3+(1-\tfrac{k}{n})^3}{\frac{k}{n}(1-\frac{k}{n})} \right) \\
\leq & 4 n \exp\left( -\frac{\varepsilon^2}{128} \cdot \beta' n  \right) ,
\end{align*}
where for the last inequality, note that $x\mapsto x^3+(1-x)^3$ is minimal on $[0,1]$ in $\tfrac{1}{2}$. The almost sure convergence of $d_K(\hat\mu_n, \mu_n)$ follows then from an application of the Borel-Cantelli Lemma.
\QED

\begin{lem}\label{gammaldp}
Let $G_1\sim \operatorname{Gamma}(\tilde\alpha_n) $ and $G_2 \sim \operatorname{Gamma}(\tilde\beta_n)$ denote independent random variables and let $\alpha_n$ and $\beta_n$ be positive real numbers going to $\infty$ such that
\begin{align*}
\lim_{n\to \infty} \frac{\tilde\alpha_n}{\alpha_n} = \alpha_0 \in (0,\infty),\qquad \lim_{n\to \infty} \frac{\tilde\beta_n}{\beta_n} = \beta_0 \in (0,\infty) .
\end{align*}
\begin{itemize}
\item[(i)] If $\beta_n/\alpha_n \to \infty$, then the vector
\begin{align*}
\mathcal{G}_n= \left( \frac{1}{\alpha_n} G_1, \frac{1}{\sqrt{\alpha_n \beta_n}} G_1, \frac{1}{\beta_n} G_1, \frac{1}{\beta_n} G_2 \right)
\end{align*}
satisfies an LDP with speed $\alpha_n$ and good rate function
\begin{align*}
I(x_1,x_2,x_3,x_4) = \alpha_0 \cdot g(\alpha_0^{-1} x_1)
\end{align*}
if $x_1>0, x_2=x_3=0, x_4 = \beta_0$ and $I(x_1,x_2,x_3,x_4)=\infty$ otherwise.
\item[(ii)] If for another set $\gamma_n$ of positive real numbers with $\gamma_n/\alpha_n \to 0$ the convergence 
\begin{align*}
\frac{\tilde\alpha_n-\tilde\beta_n}{\sqrt{\gamma_n\alpha_n}} \xrightarrow[n \rightarrow \infty ]{} 0
\end{align*}
holds, then
\begin{align*}
\mathcal{G}_n = \left( \frac{1}{\alpha_n} G_1,\frac{1}{\alpha_n} G_2, \frac{1}{2\sqrt{\gamma_n\alpha_n }} (G_1-G_2)\right)
\end{align*}
satisfies an LDP with speed $\gamma_n$ and good rate function
\begin{align*}
I(x_1,x_2,x_3) = \alpha_0^{-1} x_3^2
\end{align*}
if $x_1=x_2 = \alpha_0$ and $I(x_1,x_2,x_3)=\infty$ otherwise.
\end{itemize}
\end{lem}

\medskip

\proof In the situation (i), the logarithmic moment generating of $\alpha_n \mathcal{G}_n$ is given by 
\begin{align*}
 \log E\left[ \exp\left\{ \alpha_n \langle t,\mathcal{G}_n \rangle \right\} \right]
=&  \log E\left[ \exp\left\{ G_1\left( t_1+ t_2 \sqrt{\frac{\alpha_n}{\beta_n}} + t_3 \frac{\alpha_n}{\beta_n}\right) + G_2\cdot t_4\frac{\alpha_n}{\beta_n}  \right\} \right] \\
=&  \log\left( 1-t_1-t_2\sqrt{\frac{\alpha_n}{\beta_n}} - t_3 \frac{\alpha_n}{\beta_n}\right)^{-\tilde\alpha_n} + \frac{\tilde\beta_n}{\beta_n} \log\left(1-t_4 \frac{\alpha_n}{\beta_n}\right)^{-\beta_n} ,
\end{align*}
so we obtain the limit
\begin{align*}
\lim_{n\to \infty} \alpha_n^{-1} \log E\left[ \exp\left\{ \alpha_n \langle t,\mathcal{G}_n \rangle \right\} \right]  = -\alpha_0 \log(1-t_1) + \beta_0 t_4=: \mathcal{C}(t) .
\end{align*} 
The function $\mathcal{C}$ is finite and differentiable on $(-\infty,1)\times \mathbb{R}^3$ and steep. Then the G\"{a}rtner-Ellis Theorem (see \cite{demzei1998}) yields an LDP for $\mathcal{G}_n$ with speed $\alpha_n$ and good rate function
\begin{align*}
I(x) = \sup_{t \in \mathbb{R}^4} \left\{ \langle t,x\rangle - \mathcal{C}(t)\right\} ,
\end{align*}
which is infinite unless $x_1>0, x_2=x_3=0, x_4 = \beta_0$, in which case
\begin{align*}
I(x) = \sup_{t_1 \in \mathbb{R}} \left\{ t_1x_1 + \alpha_0 \log(1-t_1) \right\} = x_1-\alpha_0 - \alpha_0\log(\alpha_0^{-1} x_1) . 
\end{align*}
This proves part (i). For part (ii), we get for the logarithmic moment generating function
\begin{align*}
&\ \ \ \ \gamma_n^{-1} \log E\left[ \exp\left\{ \gamma_n \langle t,\mathcal{G}_n \rangle \right\} \right] \\
&=  \gamma_n^{-1} \log E\left[ \exp\left\{G_1 \left(t_1\frac{\gamma_n}{\alpha_n} + t_3 \frac{1}{2} \sqrt{\frac{\gamma_n}{\alpha_n}} \right) + G_2 \left(t_2\frac{\gamma_n}{\alpha_n} - t_3 \frac{1}{2} \sqrt{\frac{\gamma_n}{\alpha_n}} \right)\right\} \right] \\
&= \gamma_n^{-1} \log \left( 1-t_1\frac{\gamma_n}{\alpha_n} - t_3 \frac{1}{2} \sqrt{\frac{\gamma_n}{\alpha_n}}\right)^{-\tilde\alpha_n} + 
\gamma_n^{-1} \log \left( 1-t_2\frac{\gamma_n}{\alpha_n} + t_3 \frac{1}{2} \sqrt{\frac{\gamma_n}{\alpha_n}}\right)^{-\tilde\beta_n} \\
&= \gamma_n^{-1} \log \left( \left( 1-t_1\frac{\gamma_n}{\alpha_n} - t_3 \frac{1}{2} \sqrt{\frac{\gamma_n}{\alpha_n}}\right) \left( 1-t_2\frac{\gamma_n}{\alpha_n} + t_3 \frac{1}{2} \sqrt{\frac{\gamma_n}{\alpha_n}}\right) \right)^{-\tilde\alpha_n} \\
& \qquad      + \log\left( 1- t_2\frac{\gamma_n}{\alpha_n} + t_3 \frac{1}{2} \sqrt{\frac{\gamma_n}{\alpha_n}}\right)^{\gamma_n^{-1}(\tilde\alpha_n-\tilde\beta_n)} .
\end{align*}
Due to the assumption $\frac{\tilde\alpha_n-\tilde\beta_n}{\sqrt{\gamma_n\alpha_n}} \to 0$, the second term goes to zero and since $\gamma_n/\alpha_n \to 0$, the first term can be written as
\begin{align*}
  - \frac{\tilde\alpha_n}{\alpha_n} \log \left( 1-t_1\frac{\gamma_n}{\alpha_n} -t_2 \frac{\gamma_n}{\alpha_n} -t_3^2 \frac{1}{4} \frac{\gamma_n}{\alpha_n} +o\left(\frac{\gamma_n}{\alpha_n}\right) \right)^{\alpha_n/\gamma_n} ,
\end{align*}
which has the limit 
\begin{align*}
\mathcal{C}(t) := \alpha_0 (t_1+t_2+\tfrac{1}{4}t_3^2) . 
\end{align*}
As in (i), this implies an LDP for $\mathcal{G}_n$ with speed $\gamma_n$ and the stated rate function.
\QED

\medskip

\subsection{Proof of Theorem \ref{LLN2}:} 

As describes in Section 3, the weak convergence of $\mu_n$ follows, if we show entrywise convergence of the rescaled tridiagonal matrix 
\begin{align*}
(1+\sigma) \sqrt{b_n \frac{1+\sigma}{\sigma n\beta' }} \left( \mathcal{J}_n(\beta,a_n,b_n) - \frac{\sigma}{1+\sigma} I_n\right)
\end{align*}
to the tridiagonal representation $\mathcal{SC}$ of the semicircle law. By Lemma \ref{empspectral}, this is equivalent to the weak convergence of the empirical measure $\hat\mu_n$. We start with the sub-diagonal entries, whose square is given by
\begin{align*}
b_n \frac{(1+\sigma)^3}{\sigma n\beta'} p_{2k-1} (1-p_{2k-2}) p_{2k}(1-p_{2k-1} ) .
\end{align*}
The random variable $ b_n \frac{1+\sigma}{n\beta'} p_{2k}$ has mean
\begin{align*}
b_n \frac{1+\sigma}{n\beta'} \cdot \frac{(n-k)\beta'}{a_n+b_n-(2k-1)\beta'} \xrightarrow[n \rightarrow \infty ]{}1
\end{align*}
and by Lemma \ref{betaab} for $n$ sufficiently large,
\begin{align*}
P\left( \left| \frac{b_n}{n} p_{2k}- E\left[\frac{b_n}{n} p_{2k}\right] \right| >\varepsilon \right) \leq 4  \exp \left\{ -\frac{\varepsilon^2}{128} \cdot \frac{n^2}{b_n^2}\cdot \frac{\left((n-k)\beta' \right)^3+\left( a_n+b_n-(n+k-1)\beta' \right)^3}{\left((n-k)\beta' \right)\left( a_n+b_n-(n+k-1)\beta' \right)} \right\}
\end{align*}
The exponent is of order $-\varepsilon^2 n$, so the right hand side is summable. The Borel-Cantelli Lemma implies 
\begin{align} \label{proof1} 
b_n \frac{1+\sigma}{n\beta'} p_{2k} \xrightarrow[n \rightarrow \infty ]{} 1
\end{align} 
almost surely. Similarly, the random variables $p_{2k-1}$ and $p_{2k-2}$ concentrate exponentially fast around their mean, which converges to $\tfrac{\sigma}{1+\sigma}$ and 0, respectively. 
The almost sure convergence of the subdiagonal entries to 1 follows. The diagonal entries of the rescaled matrix are
\begin{align*}
 \sqrt{b_n \frac{(1+\sigma)^3}{\sigma n\beta' }} p_{2k-2}(1-p_{2k-3}) + 
 \sqrt{b_n \frac{(1+\sigma)^3}{\sigma n\beta' }} \left( p_{2k-1} - \frac{\sigma}{1+\sigma}\right)+
 \sqrt{b_n \frac{(1+\sigma)^3}{\sigma n\beta' }} p_{2k-1} p_{2k-2} .
\end{align*}
The almost sure convergence in \eqref{proof1} implies that the first and last summand vanish. For the second summand, note that 
\begin{align*}
P\left( \sqrt{\frac{b_n}{n}} \left| p_{2k-1}- E\left[p_{2k-1}\right] \right| >\varepsilon \right) \leq 4  \exp \left\{ -\frac{\varepsilon^2}{128} \cdot \frac{n }{b_n}\cdot \frac{\left(a_n-(k-1)\beta' \right)^3+\left( b_n-(k-1)\beta' \right)^3}{\left(a_n-(k-1)\beta' \right)\left( b_n-(k-1)\beta' \right)} \right\}
\end{align*}
and the right hand side is again summable. The assumption on the convergence speed of $a_n/b_n$ to $\sigma$ in Theorem \ref{LLN2} guarantees  that 
\begin{align*}
\sqrt{\frac{b_n}{n}}\left( E\left[p_{2k-1}\right] - \frac{\sigma}{1+\sigma} \right)
= \sqrt{\frac{b_n}{n}}\left( \frac{a_n-(k-1) \beta'}{a_n+b_n-(2k-2)\beta'} - \frac{\sigma}{\sigma+1}\right)
= \sqrt{\frac{b_n}{n }}\frac{(\frac{a_n}{b_n} - \sigma) + \frac{\sigma-1}{b_n}(k-1)\beta' }
{(\frac{a_n}{b_n} +1 - \frac{2k-2}{b_n}\beta' )(1+\sigma)} 
\end{align*}
vanishes as $n\to \infty$ and therefore the second summand goes to zero almost surely. Consequently all diagonal entries go to zero and the subdiagonal entries converge to 1. The entrywise limit of $(1+\sigma) \sqrt{b_n \frac{1+\sigma}{\sigma n\beta' }} \left( \mathcal{J}_n(\beta,a_n,b_n) - \frac{\sigma}{1+\sigma} I_n\right)$ is therefore the infinite matrix with spectral measure $\operatorname{SC}$. This concludes the proof of Theorem \ref{LLN1}. 
\QED

\medskip

\subsection{Proof of Theorem \ref{LLN1}:} 

The almost sure weak convergence of $\mu_n$ is a direct consequence of Theorem \ref{LDP1}, since the rate function in Theorem \ref{LDP1} is good with the unique minimizer $\operatorname{MP}(\tau)$ and the topology on $\mathcal{M}_1$ is metrizable. Therefore, $\mu_n \xrightarrow[n \rightarrow \infty ]{} \operatorname{MP}(\tau)$ and the convergence of $\hat\mu_n$ with the same limit follows from Lemma \ref{empspectral}. \QED

\medskip

\subsection{Proof of Theorem \ref{LDP2}:} 

We start by showing an LDP for the independent beta distributed entries of the rescaled matrix 
\begin{align*}
4 \sqrt{\frac{b_n}{n\beta}}\left(  \mathcal{J}_n(\beta, a_n,b_n)  - \frac{1}{2} I_n \right) .
\end{align*}
with diagonal elements
\begin{align*}
\left( 4\sqrt{\frac{b_n}{n\beta}} p_{2k-2}\right)(1-p_{2k-3}) 
+ 4 \sqrt{\frac{b_n}{n\beta}}(p_{2k-1} -\tfrac{1}{2}) 
+ p_{2k-1}\left( 4\sqrt{\frac{b_n}{n\beta}} p_{2k-2}\right)
\end{align*}
and off-diagonal elements the square root of
\begin{align*}
p_{2k-1}(1-p_{2k-2})\left(16 \frac{b_n}{n\beta}p_{2k}\right)(1-p_{2k-1}) .
\end{align*}
Recalling the definition of the random variables $p_k$ in \eqref{defbetaentries} and the representation of the beta distribution in \eqref{betafraction}, we have for the canonical moments of odd index
\begin{align*}
\left( p_{2k-1},4 \sqrt{\frac{b_n}{n\beta}} (p_{2k-1}-\tfrac{1}{2}) \right) \stackrel{d}{=} \psi \left( \frac{1}{b_n} G_1,\frac{1}{b_n} G_2,\frac{1}{2\sqrt{b_n n\beta'}}(G_1-G_2) \right) ,
\end{align*}
where $G_1\sim \operatorname{Gamma}(a_n-(k-1)\beta')$ independent of $G_2\sim \operatorname{Gamma}(b_n-(k-1)\beta')$ and 
\begin{align*}
\psi(x_1,x_2,x_3) = \left( \frac{x_1}{x_1+x_2},\frac{4}{\sqrt{2}}\frac{x_3}{x_1+x_2} \right) .
\end{align*}
The assumptions of Lemma \ref{gammaldp} (ii) are satisfies, if we set $\gamma_n=n\beta'$ and $\alpha_n=b_n$, such that $\alpha_0=1$ and 
\begin{align*}
\frac{\tilde\alpha_n-\tilde\beta_n}{\sqrt{\gamma_n\alpha_n}} = \frac{a_n-b_n}{\sqrt{n\beta' b_n}}
\xrightarrow[n \rightarrow \infty ]{} 0 .
\end{align*}
We get that $(\frac{1}{b_n} G_1,\frac{1}{b_n} G_2,\frac{1}{2\sqrt{b_n n\beta'}}(G_1-G_2))$ satisfies an LDP with speed $n\beta'$ and good rate function $I_0(x_1,x_2,x_3) = x_3^2$ if $x_1=x_2=1$ and  $I_0(x_1,x_2,x_3) = \infty$ otherwise. By the contraction principle, $(p_{2k-1},4 \sqrt{\frac{b_n}{n\beta}} (p_{2k-1}-\tfrac{1}{2}))$ satisfies the LDP with the same speed and good rate function 
\begin{align*}
I_1(y_1,y_2) = \inf \{ I_0(x_1,x_2,x_3) | \psi(x_1,x_2,x_3)=(y_1,y_2) \} = \tfrac{1}{2} y_2^2 ,
\end{align*}
if $y_1=\frac{1}{2}$ and $I_1(y_1,y_2)=\infty$ otherwise. This is immediate from the fact that on the set where $I_1$ (or $I_0$) is finite, $y_2=\tfrac{2}{\sqrt{2}}x_3$. Among the canonical moments with even index, we need to control
\begin{align*}
(y_1,y_2,y_3) := \left( p_{2k},4 \sqrt{\frac{b_n}{n\beta}} p_{2k}, 16 \frac{b_n}{n\beta}p_{2k} \right) \stackrel{d}{=} \psi \left( \frac{1}{n\beta'} G_1,\frac{1}{\sqrt{b_n n\beta'}} G_1,\frac{1}{b_n}G_1,\frac{1}{b_n}G_2 \right) ,
\end{align*}
again with independent $G_1\sim \operatorname{Gamma}((n-k)\beta')$, $G_2\sim \operatorname{Gamma}(a_n+b_n-(n-k-1)\beta')$ and 
\begin{align*}
\psi(x_1,x_2,x_3,x_4) = \left( \frac{x_3}{x_3+x_4},\frac{4}{\sqrt{2}}\frac{x_2}{x_3+x_4},8 \frac{x_1}{x_3+x_4} \right) .
\end{align*}
Lemma \ref{gammaldp}, part $(i)$, this time with $\alpha_n=n\beta', \beta_n=b_n, \alpha_0=1,\beta_0=2$, shows that the vector of rescaled gamma distributed random variables satisfies an LDP with speed $n\beta'$ and rate function 
\begin{align*}
I_2(x_1,x_2,x_3,x_4) =g(x_1)
\end{align*}
if $x_1>0,x_2=x_3=0,x_4=2$ and $I_2(x_1,x_2,x_3,x_4)=\infty$ otherwise. The contraction principle gives an LDP for $(y_1,y_2,y_3)$, where the rate function $I_3$ is finite only if $y_1=y_2=0, y_3>0$ and in this case, $y_3=4x_1$ and 
\begin{align*}
I_3(y_1,y_2,y_3) =  g\left(\frac{y_3}{4} \right) .
\end{align*}
Since $p_1,p_2,\dots $ are independent, the Dawson-G\"artner Theorem (Theorem 4.6.1 in \cite{demzei1998}) yields an LDP for the infinite vector 
\begin{align*}
y= (y_{1,1},y_{1,2},y_{2,1},y_{2,2},y_{2,3},y_{3,1},\dots ) =
 \left(p_{1},4 \sqrt{\frac{b_n}{n\beta}} (p_{1}-\tfrac{1}{2}),p_{2},4 \sqrt{\frac{b_n}{n\beta}} p_{2}, 16 \frac{b_n}{n\beta}p_{2},p_3,\dots \right)
\end{align*} 
with speed $n\beta'$ and good rate function 
\begin{align*}
I_y =  \sum_{k=1}^\infty I_1(y_{2k-1,1},y_{2k-1,2}) + I_3(y_{2k,1},y_{2k,2},y_{2k,3}) ,  
\end{align*}
finite only if $y_{2k-1,1}=\tfrac{1}{2},y_{2k,1}=y_{2k,2}=0$ and $y_{2k,3}>0$ and in this case,
\begin{align*}
I_y(y) = \sum_{k=1}^\infty \frac{1}{2} y_{2k-1,2}^2 +  g\left(\frac{y_{2k,3}}{4} \right) .
\end{align*} 
The entries of the rescaled tridiagonal matrix, which are the recursion coefficients of the spectral measure, are given by
\begin{align*}
c_k^2 =& y_{2k-1,1}(1-y_{2k-2,1}) y_{2k,3}(1-y_{2k-1,1}), \\
d_k =& y_{2k-2,2}(1-y_{2k-3,1}) + y_{2k-2,2}y_{2k-1,1} + y_{2k-1,2} .
\end{align*}
In particular, the recursion coefficients are a continuous function of $y$. Applying the contraction principle once more, we can set $c_k^2=\tfrac{1}{4} y_{2k,3}$ and $d_k=y_{2k-1,2}$ and obtain that
the sequence $r=(d_1,c_1,d_2,\dots )$ of recursion coefficients satisfies the LDP with speed $n\beta'$ and good rate
\begin{align*}
I_r(d_1,c_1,d_2,\dots )= \sum_{k=1}^\infty \frac{1}{2} d_k^2 +  g\left( c_k^2 \right) .
\end{align*} 
To transfer this LDP to the spectral measure, let $\mathcal{R}_c$ denote the set of sequences $r=r(\mu)$ of measures $\mu \in \mathcal{M}_c$ with compact support. The mapping $r(\mu)\mapsto \mu$ is well-defined on $\mathcal{R}_c$ and continuous and so the contraction principle yields the LDP for the spectral measure $\mu_n$ with rate function $\mathcal{I}(\mu) = I_r(r(\mu)) = \mathcal{I}_L(\mu)$.
\QED

\medskip

\subsection{Proof of Theorem \ref{LDP1}:}

We need to prove an LDP for the entries of the rescaled matrix
\begin{align*}
\frac{b_n}{a_n} \mathcal{J}_n(\beta, a_n,b_n) ,
\end{align*}
which has diagonal entries
\begin{align*}
\frac{b_n}{a_n} p_{2k-2} (1-p_{2k-3}) + \frac{b_n}{a_n} p_{2k-1}(1-p_{2k-2})
\end{align*}
and the squared off-diagonal entries are
\begin{align*}
 \frac{b_n}{a_n} p_{2k-1}(1-p_{2k-2})\cdot \frac{b_n}{a_n} p_{2k}(1-p_{2k-1}) .
\end{align*}
The canonical moments with odd index appear in the matrix in two ways, and can be written as
\begin{align*}
\left( p_{2k-1},\frac{b_n}{a_n}p_{2k-1} \right) \stackrel{d}{=} \phi \left( \frac{1}{a_n} G_1,\frac{1}{b_n} G_1,\frac{1}{b_n}G_2 \right) ,
\end{align*}
where $G_1\sim \operatorname{Gamma}(a_n-(k-1)\beta')$ independent of $G_2\sim \operatorname{Gamma}(b_n-(k-1)\beta')$ and 
\begin{align*}
\phi(x_1,x_2,x_3) = \left( \frac{x_2}{x_2+x_3},\frac{x_1}{x_2+x_3} \right) .
\end{align*}
From the first part of Lemma \ref{gammaldp} with $\alpha_0=\beta_0=1$, we see that the random vector $ (\frac{1}{a_n} G_1,\frac{1}{b_n} G_1,\frac{1}{b_n}G_2)$ satisfies an LDP
with speed $a_n$ and rate function $I_0 (x_1,x_2,x_3)= g(x_1)$ if $x_1>0, x_2=0,x_3=1$ and $I_0 (x_1,x_2,x_3)=\infty$ otherwise. Then the contraction principle gives an LDP for $(  p_{2k-1},\frac{b_n}{a_n}p_{2k-1})$ with rate function
\begin{align*}
I_1(y_1,y_2) = \inf \{I_0(x_1,x_2,x_3) | \phi(x_1,x_2,x_3) = (y_1,y_2)\} = g(y_2)
\end{align*}
if $y_2>0$ and $y_1=0$ and $I_1(y_1,y_2)=\infty $ otherwise. Turning to the random variables with even index, we have that
\begin{align*}
\left( p_{2k},\frac{b_n}{a_n}p_{2k} \right) \stackrel{d}{=} \phi \left( \frac{1}{a_n} G_1,\frac{1}{b_n} G_1,\frac{1}{b_n}G_2 \right) ,
\end{align*}
with $G_1\sim \operatorname{Gamma}((n-k)\beta')$ and $G_2\sim \operatorname{Gamma}(a_n+b_n-(n-k-1)\beta')$ independent and $\phi$ as above. Applying again the first part of Lemma \ref{gammaldp} with $\alpha_0=\tau$ and $\beta_0$, we get an LDP for the vector of gamma distributed random variables with speed $a_n$ and good rate $I_2(x_1,x_2,x_3)= \tau g(x_1/\tau)$ for $x_1>0,x_2=0,x_3=1$. The contraction principle gives the LDP for $(  p_{2k},\frac{b_n}{a_n}p_{2k})$ with rate function $I_3(y_1,y_2) = \tau I_1(y_1,y_2/\tau )$. 
Collecting the independent canonical moments we get a projective LDP for
\begin{align*}
y= (y_{1,1},y_{1,2},y_{2,1},y_{2,2},\dots ) = \left( p_{1},\frac{b_n}{a_n}p_{1},p_{2},\frac{b_n}{a_n}p_{2},\dots  \right)
\end{align*}
with speed $a_n$ and good rate function
\begin{align*}
I_y(y) = \sum_{k=1}^\infty I_1(y_{2k-1,1},y_{2k-1,2})+\tau I_1(y_{2k,1},y_{2k,2}/\tau ) .
\end{align*}
We now perform a transformation to the random variables $z_{k}=y_{k,2}(1-y_{k,1})=\frac{b_n}{a_n}p_k(1-p_{k-1})$ with $p_0=0$. On the set where the rate function $I_y$ is finite we have $y_{k,1}=0$ and $z_k=y_{k,2}$ such that $z=(z_1,z_2,\dots)$ satisfies the LDP with rate
\begin{align*}
I_z(z) = \sum_{k=1}^\infty I_1(0,z_{2k-1}) + \tau I_1(0,z_ {2k}/\tau ) .
\end{align*}
To complete the proof, note that $z$ contains the recursion variables of the spectral measure $\mu_n$ of the rescaled random matrix. It remains to apply the continuous mapping $z(\mu)\mapsto \mu$ from the recursion variables to the spectral measure analogous to the step in Section 4.3 to complete the proof.
\QED

\medskip

\bibliographystyle{apalike}
\bibliography{detnag}

\end{document}